\documentclass{article}

\usepackage{Jared-MK}
\usepackage{graphicx} 
\usepackage{mathtools} 
\usepackage{hyperref}
\usepackage{thmtools}
\usepackage{thm-restate}
\usepackage{cleveref}

\newtheorem{assumption}{Assumption}
\newcommand{\eps}{\epsilon}

\newcommand{\p}{\partial}
\newcommand{\n}{\nabla}
\newcommand{\tX}{\tilde{X}}
\newcommand{\EGLe}{E_{\text{GL}, \eps}}
\newcommand{\EAYe}{E_{\text{AYMH}, \eps}}
\newcommand{\eGLe}{e_{\text{GL}, \eps}}
\newcommand{\eAYe}{e_{\text{AYMH}, \eps}}

\newtheorem*{theorem*}{Theorem}

\newcommand{\sauce}{\textcolor{red}{source here}}

\title{Second Inner Variations of Energy and Index of Codimension $2$ Minimal Submanifolds}

\author{Jared Marx-Kuo}
\date{July 18, 2023}
\begin{document}

\maketitle

\begin{abstract}
\noindent We compute the second inner variation of the Abelian Yang--Mills--Higgs and Ginzburg--Landau energies. Given a sequence of critical points with energy measures converging to a codimension $2$ minimal submanifold, we use the second inner variation formula to bound the morse index of the submanifold by the index of the critical points. The key tools are the convergence of the energy measures and the stress-energy tensors of solutions to Abelian Yang--Mills--Higgs and Ginzburg--Landau equations.
\end{abstract}
\tableofcontents

\section{Introduction}
There is a long history of approximating minimal submanifolds with solutions to partial differential equations. The Allen--Cahn energy gives one such equation
\begin{align} \label{ACEnergy}
E_{AC,\eps}(u) &= \int_M \frac{\eps}{2} |\nabla u|^2 + \frac{W(u)}{\eps} \\ \label{ACEquation}
\eps^2 \Delta_g(u) &= u(u^2 -1) \iff E_{AC,\eps}'(u) = 0
\end{align}
Modica--Mortola first established $\Gamma$-convergence of the Allen--Cahn energy to the perimeter functional \cite{mortola1977esempio} and Hutchison--Tonegawa later showed weak convergence of the energy measures to a stationary codimension $1$ varifold \cite{hutchinson2000convergence}. Pacard--Ritor\`e \cite{pacard2003constant} then showed that given a non-degenerate minimal surface, an Allen--Cahn solution can be glued in nearby. Guaraco established a min-max construction for the Allen--Cahn energy \cite{guaraco2018min}, and Chodosh--Mantoulidis established curvature estimates and index bounds using the Allen--Cahn equation \cite{chodosh2020minimal}. \nl \nl
These results highlight a deep connection between Allen--Cahn solutions and minimal surfaces. In particular, we are interested in the following: for $\{u_{\eps}\}$ a sequence of critical points to equation \eqref{ACEnergy}, Le \cite{le2011second} and Gaspar \cite{gaspar2020second} computed the second inner variation of the $E_{AC,\eps}(u)$ and use it to bound the area Morse index of a limiting minimal surface, $Y^{n-1} \subseteq M^n$, by the Allen--Cahn energy morse index of the solutions $\{u_{\eps}\}$. The same result was also proved by Hiesmayr \cite{hiesmayr2018spectrum} but using different techniques:
\begin{theorem}[Gaspar 2018, Hiesmayr 2018] \label{ACSecondInner}
Let $M$ a closed Riemannian manifold of dimension $n \geq 3$ and $\{u_{\eps_k}\}$ a sequence of critical points of equation \eqref{ACEnergy} with $\eps_k \downarrow 0$. Assume that there are positive constants $c_0$, $E_0$, and a non-negative integer $p$ such that 
\[
\limsup_k \sup_M |u_{\eps_k}| \leq c_0, \qquad \limsup_k E_{\eps_k}(u_{\eps_k}) \leq E_0, \qquad \limsup_k m(u_{\eps_k}) \leq p
\]
Then up to a subsequence, for all open subsets $U \subset \subset M \backslash \text{sing}(\Gamma)$ with $U \cap \text{reg}(\Gamma) \neq \emptyset$, the eigenvalues $\{\lambda_{\ell}^{\eps}(U)\}_{\ell}$ of the linearized Allen--Cahn operator at $u_{\eps_l}\Big|_U$, $L_k = - \eps_k \Delta + W''(u_{\eps_k})/\eps_k$, and the eigenvalues $\{\lambda_{\ell}(U)\}_{\ell}$ of the Jacobi operator of $\Gamma$ acting on normal vector fields supported on $U \cap \Gamma$ satisfy
\[
\limsup_{k} \frac{\lambda_{\ell}^{\eps_k}(U)}{\eps_k} \leq \lambda_{\ell}(U)
\]
for all $\ell$. In particular, $\text{reg}(\Gamma)$ has morse index at most $p$.
\end{theorem}
\noindent Many of the results connecting the Allen--Cahn equation to minimal hypersurfaces are now being extended to the codimension $2$ case. In this context, the Abelian Yang--Mills--Higgs equations (sometimes referred to as ``$U(1)$-Yang--Mills--Higgs") are
\begin{align} \label{AYMHEquation}
\eps^2 \nabla^* \nabla u &= \frac{1}{2} (1 - |u|^2)u \\ \nonumber
\eps^2 d^* F_{\nabla} &= i \text{Im}\langle u, \nabla u \rangle_L 
\end{align}
and Ginzburg--Landau equations 
\begin{equation} \label{GLEquation}
\eps^2 \Delta u = u(|u|^2 - 1)
\end{equation}
These equations have been used to approximate minimal submanifolds $Y^{n-2}$. In the Abelian Yang--Mills--Higgs setting, Jaffe and Taubes \cite{jaffe1980monopoles} initially constructed solutions to equation \eqref{AYMHEquation} when $M = \R^2$, concentrating at a set of points  $\{z_i\}_{i = 1}^k$ as $\eps \to 0$. In all dimensions, Badran--Del-Pino \cite{badran2022entire} \cite{badran2022solutions} have completed a gluing construction, paralleling the result of Pacard--Ritore. Parise--Pigati--Stern \cite{parise2021convergence} have established $\Gamma$-convergence properties and Pigati--Stern \cite{pigati2021minimal} established weak convergence of the energy measures to stationary \textit{integral} rectifiable $(n-2)$-varifolds:
\begin{theorem}[Pigati-Stern, Thm 1.1] \label{PSAYMHTheorem}
Let $\{(u_{\eps}, \nabla_{\eps})\}$ a family of critical points of solutions to equation \eqref{AYMHEquation} satisfying
\[
E_{\eps}(u_{\eps}, \nabla_{\eps}) \leq \Lambda < \infty
\]
Then as $\eps \to 0$, the energy measures
\[
\mu_{\eps} := \frac{1}{2\pi} e_{\eps}(u_{\eps}, \nabla_{\eps}) vol_g
\]
converge subsequentially, in duality with $C^0(M)$, to the weight measure of a stationary integral $(n-2)$-varifold $V$. Also for all $0 \leq \delta < 1$, 
\[
\text{spt}(V) = \lim_{\eps \to 0} \{ |u_{\eps}| \leq \delta \}
\]
in the Hausdorff topology.
\end{theorem}
\noindent In the Ginzburg--Landau setting, Lin--Riviere \cite{lin2001quantization},  Betheul--Brezis--Orlandi \cite{bethuel2001asymptotics}, and Betheul--Orlandi--Smets \cite{bethuel2006convergence} showed that solutions concentrate about a collection of \textit{codimension $2$ minimal submanifolds}. This was later refined in the Riemannian setting by Cheng \cite{ChengAsymptotics}, Stern \cite{stern2021existence} and Pigati--Stern \cite{pigati2022quantization}. These authors show that the energy of a sequence of solutions converge to a stationary, but potentially \textit{not integral}, $n-2$ varifold, plus a diffuse measure on all of $M$. 
\begin{theorem}[Betheul--Brezis--Orlandi, Cheng, Stern, Pigati--Stern, Betheul--Orlandi--Smets] \label{PSGLTheorem}
Let $(M^n, g)$ a closed manifold with $n \geq 3$ and $\{u_{\eps}\}$ a sequence of solutions to equation \eqref{GLEquation} with 
\[
\limsup_{\eps \to 0} \frac{1}{|\log \eps|} \int_K \left(\frac{1}{2} |d u_{\eps}|_g^2 + \frac{W(u_{\eps})}{\eps^2} \right) d vol_{g_{\eps}} < \infty
\]
for all compact $K \subseteq M$. Then, up to a subsequence, the normalized energy densities
\[
\mu_{\eps} = \frac{1}{|\log \eps|} \left(\frac{1}{2} |d u_{\eps}|_g^2 + \frac{W(u_{\eps})}{\eps^2} \right)
\]
converge to a radon measure, $\mu$, which decomposes as 
\[
\mu = |V| + f vol_{g_0}
\]
for $f$ non-negative.
\end{theorem}
\noindent While the diffuse measure and lack of integrality is differs from the Abelian Yang--Mills--Higgs setting, Cheng \cite{cheng2017geometric} and Stern \cite{stern2021existence} have shown that there exists solutions, $\{u_{\eps_k}\}$, for which $f = 0$ and $V$ is integral. When the minimal submanifold is prescribed, De-Phillipis--Pigati \cite{de2022non} used variational methods to prove the existence of solutions accumulating along codimension $1$ and $2$ in each of the Allen-Cahn, Abelian Yang--Mills--Higgs, and Ginzburg--Landau settings. \nl \nl
\noindent Given the parallels between the Allen-Cahn, Abelian Yang--Mills--Higgs, and Ginzburg--Landau equations in $\Gamma$-convergence and convergence of the energy measures, one expects to replicate results which pass geometric information between a sequence of solutions, $\{(u_{\eps}, \n_{\eps})$ ( $(\{u_{\eps}\}$), and the limiting stationary varifold, $V$. This is the goal of this paper: to recreate the Morse-Index bound of Le and Gaspar in the context of the Ginzburg--Landau and Abelian Yang--Mills--Higgs equations on Riemannian manifolds. 

\section{Statement of Results}
Let $\Gamma$ be a minimal codimension $2$ submanifold with components $\{\Gamma_i\}$. Define
\[
\mu_V := \sum_{i = 1}^q m_i d \mathcal{H}^{n-2}_{\Gamma_i} \\
\]
where $m_i$ are positive integers so that $V = v(\Gamma, \mu_V)$ defines a stationary integral $(n-2)$-varifold. We will assume the following:
\begin{assumption} \label{ConvergenceAssumption}
A sequence of solutions, $\{(u_{\eps_k}, \n_{\eps_k})\}$ or $\{u_{\eps_k}\}$, to the Abelian Yang--Mills--Higgs or Ginzburg--Landau equations has bounded total energy, i.e.
\begin{align*}
\limsup_{k \to \infty} E_{AYMH, \eps_k}(u_{\eps_k}, \n_{\eps_k}) &< \infty \\
\limsup_{k \to \infty} E_{GL, \eps_k}(u_{\eps_k}) &< \infty
\end{align*} 
and converges to $V$ in that:
\begin{align} \label{AYMHConvergenceAssumption}
\frac{1}{2\pi}\lim_{k \to \infty} e_{AYMH, \eps_k}(u_{\eps_k}, \n_{\eps_k}) d Vol_M &= \mu_V \\ \label{GLConvergenceAssumption}
\frac{1}{\pi}\lim_{k \to \infty} \frac{1}{|\log(\eps_k)|} e_{GL, \eps_k}(u_{\eps_k}) d Vol_M &= \mu_V
\end{align}
\end{assumption}
\noindent Regularity aside, theorem \ref{PSAYMHTheorem} says that equation \eqref{AYMHConvergenceAssumption} is typical. By contrast, theorem \ref{PSGLTheorem} says that equation \eqref{GLConvergenceAssumption} is atypical, but Cheng \cite{cheng2017geometric} and Stern \cite{stern2021existence} show that $b_1(M) = 0$ implies $f = 0$ in theorem \ref{PSGLTheorem}. In any case, we assume both convergence assumptions to establish morse index bounds. We state our main theorems here:
\begin{theorem} \label{SecondInnerVarThm}
Let $\{(u_{\eps_k}, \n_{\eps_k})\}$ (resp. $\{u_{\eps_k}\}$) be a sequence of critical points for Abelian Yang--Mills--Higgs (resp. Ginzburg--Landau) satisfying assumption \eqref{AYMHConvergenceAssumption} (resp. \eqref{GLConvergenceAssumption}). For $X \in \Gamma(TM)$, let $\Phi^t: M \to M$ be the corresponding flow, and define $\{(u_{\eps_k}^t, \n_{\eps_k}^t)\}$ (resp. $\{(u_{\eps}^t)\}$). Then
\begin{align} \label{SecondInnerVarEqn}
\frac{1}{2\pi}\lim_{k \to \infty} \frac{d^2}{dt^2} \EAYe(u_{\eps_k}^t, \n_{\eps_k}^t) \Big|_{t = 0} &= \sum_{i = 1}^q m_i \left[ D^2A(\Gamma_i)(X,X) + \int_{\Gamma_i} \left(- \div_{N \Gamma_i}(X)^2 + \frac{1}{2} ||\dot{g}_X||_{N \Gamma_i}^2 d Vol_{\Gamma_i}\right) \right] \\ \nonumber
\frac{1}{\pi}\lim_{k \to \infty} \frac{d^2}{dt^2} \EGLe(u_{\eps_k}^t) \Big|_{t = 0} &= \sum_{i = 1}^q m_i \left[ D^2A(\Gamma_i)(X,X) + \int_{\Gamma_i} \left(- \div_{N \Gamma_i}(X)^2 + \frac{1}{2} ||\dot{g}_X||_{N \Gamma_i}^2 d Vol_{\Gamma_i}\right) \right]
\end{align}
Moreover, the integral term in \eqref{SecondInnerVarEqn} is non-negative.
\end{theorem}
\noindent \rmk \; We note that Le (\cite{le2015second}, theorem 1.5) and Cheng (\cite{ChengInstability}, Proposition 2.6) have done this computation for a sequence of Ginzburg--Landau solutions in the Euclidean and Riemannian setting, respectively. Both authors compute this to investigate stability of Ginzburg-Landau solutions. We also note that Cheng (\cite{ChengInstability}, Proposition 5.2) computes the second inner variation for Abelian Yang-Mills-Higgs in a different form, again to prove stability. We include our own derivations since they emphasize the parallels between the Abelian Yang--Mills--Higgs and Ginzburg--Landau cases, unified by the convergence of the stress energy tensor. \nl \nl
\noindent For any $U \subseteq M$ with $U \cap \text{reg}(V) \neq \emptyset$, let $\{\lambda_p(U)\}$ denote the eigenvalues of the Jacobi operator of $V$ acting on normal vector fields supported on $U \cap \text{reg}(V)$. Similarly, for a sequence of critical points $\{(u_{\eps_k}, \n_{\eps_k})\}$ (resp. $\{u_{\eps_k}\}$ ) to the Abelian Yang--Mills--Higgs (resp. Ginzburg--Landau) equations, let $\lambda_{AYMH, p}^{\eps_k}(U)$ (resp. $\lambda_{GL, p}^{\eps_k}$) denote the eigenvalues of the linearized Abelian Yang--Mill--Higgs (resp. Ginzburg--Landau) operator at $(u_{\eps_k}, \n_{\eps_k}) \Big|_U$ (resp. $u_{\eps_k} \Big|_U$).
\begin{theorem} \label{MorseIndexTheorem}
Let $\{(u_{\eps_k}, \n_{\eps_k})\}$ (resp. $\{u_{\eps_k}\}$) be a sequence of critical points for Abelian Yang--Mills--Higgs (resp. Ginzburg--Landau) satisfying assumption \ref{ConvergenceAssumption} and 
\begin{align*}
\limsup_k m(u_{\eps_k}, \n_{\eps_k})  &\leq m \\
\limsup_k m(u_{\eps_k}, \n_{\eps_k})  &\leq m 
\end{align*}
Then up to a subsequence, for all $U \subseteq M \backslash \text{sing}(V)$ with $U \cap \text{reg}(V) \neq \emptyset$
\begin{align*}
\limsup_k \lambda_{AYMH, p}^{\eps_k}(U) &\leq 2\lambda_p(U) \\
\limsup_k \lambda_{GL, p}^{\eps_k}(U) &\leq \lambda_p(U)
\end{align*}
In particular, the morse index of the regular part of $V$ has morse index at most $m$.
\end{theorem}
\noindent We can immediately apply this to any of the solutions constructed in \cite{de2022non}, \cite{cheng2017geometric}, \cite{stern2021existence}, \cite{badran2022solutions} to give a lower bound on the morse index of the constructed solutions, $\{u_{\eps}\}$ or $\{(u_{\eps}, \n_{\eps})\}$, in terms of the morse index of the limiting varifold. \nl \nl
\noindent The techniques to prove theorems \ref{SecondInnerVarThm} \ref{MorseIndexTheorem} are almost identical to that of Gaspar \cite{gaspar2020second} and Hiesmayr \cite{hiesmayr2018spectrum}, expect we replace the Allen-Cahn equation with the Abelian Yang--Mills--Higgs and Ginzburg--Landau equations. We also substitute proposition 2.2 of Gaspar \cite{gaspar2020second} with knowledge of the limit of the stress-energy tensor for the Abelian Yang--Mills-Higgs and Ginzburg--Landau equations. We also remark that up to the appropriate constant, the integral term in \eqref{SecondInnerVarEqn} simplifies to that of Gaspar's (\cite{gaspar2020second}, Proposition 3.3) and Le's (\cite{le2011second}, Theorem 1) when the normal bundle is $1$-dimensional and $X$ lies in the normal bundle. 
\subsection{Acknowledgements}
The author would like to thank Otis Chodosh for suggesting the initial idea of the project. The author would also like to thank Daniel Stern and Pedro Gaspar for their time in answering questions about their relevant papers. The author dedicates this paper to his grandmother, Shirley Kuo.

\section{Preliminaries}
\subsection{Riemannian Background}
Let $(M, g)$ a Riemannian manifold. For $X \in \Gamma(TM)$, let $\Phi^t$ be the associated flow and define $g^t = (\Phi^t)^*(g)$. Fix a point $q \in M$ and a basis $\{\p_i\}_{i = 1}^n$ at $q$. Let $\dot{g}_{ij}, \ddot{g}_{ij}$ denote the first and second derivatives for the metric coefficients evaluated at $p$, i.e.
\[
\dot{g}_{ij} = \frac{d}{dt} g^t_{ij} \Big|_{t = 0}, \qquad \ddot{g}_{ij} = \frac{d^2}{dt^2} g^t_{ij} \Big|_{t = 0}
\]
Adopt the same notation for $\dot{g}^{ij}, \ddot{g}^{ij}$ as derivatives of the metric inverse coefficients. We also define 
\begin{align} \label{TTensorEqn}
T_{P}(V,W) &:= (g|_{P})^{ij} (g|_{P})^{k \ell} \langle \n_i V, \p_k \rangle \langle \n_{\ell} W, \p_{j} \rangle
\end{align}
where $P$ is any sub-bundle of $TM$ and $g|_{P}$ is the restricted metric. Then we have:
\begin{lemma}
The following hold
\begin{align} \label{MetricDerivatives}
\dot{g}_{ij} &= \langle \n_i X, \p_j \rangle + \langle \p_i, \n_j X \rangle \\ \nonumber
\dot{g}^{ij} &= - g^{ik} \dot{g}_{k m} g^{m j} \\ \nonumber
\ddot{g}_{ij} &= \langle \n_i \n_X X, \p_j \rangle + \langle \n_j \n_X X, \p_i \rangle + 2 \langle \n_i X, \n_j X \rangle - 2 R(X, \p_i, X, \p_j) \\ \nonumber
\ddot{g}^{ij} &= 2 g^{ir} \dot{g}_{rs} g^{sk} \dot{g}_{k m} g^{m j} - g^{ik} \ddot{g}_{km} g^{mj} \\ \nonumber
\frac{d}{dt} \sqrt{\det g^t} \Big|_{t = 0}&= \div(X) \sqrt{\det g} \\ \nonumber
\frac{d^2}{dt^2} \sqrt{\det g^t} \Big|_{t = 0} &= \div(\n_X X) - \Ric(X, X) - T_{TM}(X,X) + \div(X)^2
\end{align}
\end{lemma}
\noindent See \cite{gaspar2020second}, lemma 3.1 for reference. Here we follow the convention 
\[
R(V,W,Y,Z) = \langle (\nabla_{W} \nabla_V - \n_V \n_W - \n_{[W,V]}) Y, Z \rangle
\]
We also recall that for $X$, a not necessarily normal vector field, that the second variation of area along a submanifold $\Sigma \subseteq M$, flowed by $\Phi^t$ is given by 
\begin{equation} \label{SecondVarArea}
D^2 A\Big|_{\Sigma}(X, X) = \int_{\Sigma} [ -\tr_{T \Sigma} R(X, \cdot, X, \cdot) - T_{T \Sigma}(X,X) + \div_{T \Sigma}(X)^2 + \div_{T \Sigma}(\n_X X) + ||\nabla_{T \Sigma}^{N \Sigma} X||^2 ] d Vol_{\Sigma}
\end{equation}
(See \cite{simon1983lectures}, 9.4 or \cite{colding2011course}, section 8.1).

\subsection{Ginzburg--Landau and Abelian Yang--Mills--Higgs Equations}
The Ginzburg--Landau functional is given by
\begin{align} \nonumber
u& \in H^1(M, \C)\\ \nonumber
e_{GL, \eps}(u) &:= \frac{1}{2} |\nabla u|^2 + \frac{W(u)}{\eps^2} \\ \label{GLEnergy}
E_{GL,\eps}(u) &= \int_M \frac{1}{2 } |\nabla u|^2 + \frac{W(u)}{\eps^2}  \\ \nonumber
\mu_{GL,\eps_k} & := \frac{1}{\pi |\log(\eps)|} e_{\eps_k}(u_{\eps_k}) d Vol_g
\end{align}
with critical points of \eqref{GLEnergy} satisfying equation \eqref{GLEquation}. 
\[
\eps^2 \Delta u = u(|u|^2 - 1)
\]
%
%
\noindent For Abelian Yang--Mills--Higgs, let $L$ be a complex line bundle over $M$, then
\begin{align} \nonumber
u&: M \to L \\ \nonumber
\nabla&: \Gamma(TM \otimes L) \to \Gamma(L) \\ \nonumber
e_{AYMH, \eps}(u, \nabla)&:= |\nabla u|^2 + \eps^2 ||F_{\nabla}||^2 + \frac{W(u)}{\eps^2} \\ \label{AYMHEnergy}
\EAYe(u, \nabla)&:= \int_M |\nabla u|^2 + \eps^2 ||F_{\nabla}||^2 + \frac{W(u)}{\eps^2} \\ \nonumber
\mu_{AYMH,\eps_k} & := \frac{1}{2\pi} e_{\eps_k}(u_{\eps_k}, \n_{\eps_k}) d Vol_g
\end{align}
Throughout this paper, we will refer to $\mu_{GL, \eps_k}$ and $\mu_{AYMH, \eps_k}$ as the ``energy measures." Here, we follow the weighting and norm convention for 2-forms as in \cite{pigati2021minimal}: let $\{x_i\}$ be a local coordinate basis, then
\begin{align*}
	\tau &= \sum_{i < j} \tau_{ij} dx_i \wedge dx_j \\
	||\tau||_g^2 &= g^{ij} g^{kl} \tau_{ik} \tau_{j \ell}
\end{align*}
when $\{\p_{i}\}$ are orthonormal at a point, this gives
\[
||\tau||_g^2 = \sum_{i < j} \tau_{ij}^2 = \frac{1}{2} \sum_{i, j = 1}^n \tau_{ij}^2
\]
We also define an $\eps$-weighted inner product on pairs of sections and $i \R$-valued one forms
\begin{equation} \label{EpsInnerProduct}
\langle (f, a), (h, b) \rangle := \langle f, h \rangle_L + \eps^2 \langle a, b \rangle_L
\end{equation}
Critical points of equation \eqref{AYMHEnergy} satisfy the coupled system equation \eqref{AYMHEquation}
\begin{align*}
\eps^2 \nabla^* \nabla u &= \frac{1}{2} (1 - |u|^2)u \\
\eps^2 d^* F_{\nabla} &= i\text{Im}\left(\langle u, \nabla u \rangle_L \right)
\end{align*}
\subsubsection{Stress Energy Tensors}
For each equation, we define the corresponding stress and stress-energy tensors, following the convention of \cite{pigati2021minimal}, section $4$.
\begin{align} \label{GLTensors}
S_{GL, \eps} &= (\nabla u_{\eps}^* \nabla u) \\ \nonumber
S_{GL, \eps}(V, W) &:= \langle \nabla_V u_{\eps}, \nabla_W u_{\eps} \rangle_{\C} \\ \nonumber
T_{GL, \eps} &:= e_{\eps}(u) g - S_{GL, \eps}
\end{align}
For Abelian Yang--Mills--Higgs, the gauge group of $L$ is $U(1)$. We can write $F_{\n} = i \omega$ for $\omega \in \Omega^2(TM, \R)$. We then define
\begin{align*} \label{AYMHTensors}
(\nabla u^* \nabla u)(V,W) &:= \langle \nabla_{V} u, \nabla_{W} u \rangle_L \\
(\omega^* \omega)(V,W) &:= g^{ij} \omega(V, e_i) \omega(W, e_i) \\
S_{AYMH, \eps} &:= 2 [\nabla u^* \nabla u + \omega^* \omega] \\
T_{AYMH, \eps} &:= e_{\eps}(u,A) g - S_{AYMH, \eps}
\end{align*}
From Pigati--Stern \cite{pigati2021minimal}, we recall the convergence of the stress energy tensor
\begin{proposition}[Pigati--Stern, Proposition 6.4] \label{AYMHStressEnergyTensorConvergence}
For a family $\{(u_{\eps}, \n_{\eps})\}$ of solutions to equation \eqref{AYMHEquation} with uniform energy bound, after passing to a subsequence $\{\eps_j\} \to 0$, there exists a stationary, rectifiable, integral $(n-2)$-varifold, $V = v(\Sigma^{n-2}, \theta)$ such that
\[
\lim_{\eps \to 0} \int_M \langle T_{AYMH,\eps}(u_{\eps}, \n_{\eps}), P \rangle = \int_{\Sigma} \theta(x) \langle T_x \Sigma, P(x) \rangle d \mathcal{H}^{n-2}
\]
for every $P \in C^0(M, \text{Sym}(TM))$. 
\end{proposition}
%
\noindent As a result of proposition \ref{AYMHStressEnergyTensorConvergence} and theorem \ref{PSAYMHTheorem}, we see that for a sequence $\{(u_{\eps_k}, \n_{\eps_k})\}$ of solutions to equation \eqref{AYMHEquation}, we have
\begin{equation} \label{AYMHStressTensorEq}
\lim_{k \to \infty} \langle S_{AYMH, \eps}(u_{\eps_k}, \n_{\eps_k}), P \rangle = \int_{\Sigma} \theta(x) \langle N \Sigma_x, P \rangle d \mathcal{H}^{n-2}
\end{equation}
In the Ginzburg--Landau setting, we have the following from \cite{bethuel2001asymptotics}, \S IX, and Cheng \cite{ChengInstability} (in the proof of proposition 2.6):
%
\begin{proposition} \label{GLStressEnergyTensorConvergence}
For a sequence $\{u_{\eps_k}\}$ with uniformily bounded energy converging to a stationary varifold (i.e. $h_0 = 0$ in the context of theorem \ref{PSGLTheorem}),
\begin{align} \label{GLStressEnergyTensorEq}
\lim_{k \to \infty} \int_M \langle T_{GL, \eps_k}, P \rangle &= \int_{\Sigma} \theta(x) \langle T_x \Sigma, P(x) \rangle d \mathcal{H}^{n-2}_{\Sigma} \\ \label{GLStressTensorEq}
\lim_{k \to \infty} \int_M \langle S_{GL, \eps_k}, P \rangle &= \int_{\Sigma} \theta(x) \langle N_x \Sigma, P(x) \rangle d \mathcal{H}^{n-2}_{\Sigma}
\end{align}
for all $P \in C^0(M, \text{Sym}(TM))$.
\end{proposition}
%
%
\noindent \rmk \; In the context of our assumptions \ref{ConvergenceAssumption}, we note that the density $\theta(x)$ will be $2 \pi \Z$ (resp. $\pi\Z$) valued in accordance with the normalization in equation \eqref{AYMHConvergenceAssumption} (resp. \eqref{GLConvergenceAssumption}).

\section{Computation of first inner variation}
In the Ginzburg--Landau setting, define
\[
u^t = (\Phi^{-t})^*(u)
\]
and we have
\begin{align*}
E_{GL, \eps}(u^t) &= \int_M \frac{1}{2} |\nabla u^t|^2 + \frac{W(u^t)}{\eps^2} dV \\
&= \int_M \left[\frac{1}{2} |\nabla u^t|^2 \circ \Phi^t + \frac{W(u)}{\eps^2} \right](\Phi^t)^*(d Vol) 
\end{align*}
Closely mimicking Gaspar \cite{gaspar2020second}, \S 3, let $\{e_i^t\}$ be an ONB at $y$. We compute 
\begin{align*}
|\nabla u^t|^2 \Big|_{y} & = \sum_{i = 1}^n d u^t( e_i^t)^2 \\
&= \sum_{i = 1}^n \left( du \Big|_{\Phi^{-t}(y)} \right)(\Phi^{-t}_*(e_i^t))^2 \\
&= \sum_{i = 1}^n \left( du \right)(v_i^t)^2 \Big|_{\Phi^{-t}(y)} \\
&= \sum_{i = 1}^n g^t\left(\nabla^{g^t} u, v_i^t\right)^2 \Big|_{\Phi^{-t}(y)} \\
\implies |\nabla u^t|^2 \circ \Phi^t (x) &= \sum_{i = 1}^n g^t\left(\nabla^{g^t} u, v_i^t\right)^2 \Big|_{x}
\end{align*}
Here, we've defined $v_i^t = \Phi^{-t}_* \Big|_{y}(e_i^t)$, which is a vector at $\Phi^{-t}(y)$. Moreover, we've noted that $v_i^t$ is an ONB at $\Phi^{-t}(y)$ with respect to the metric $g^t = (\Phi^t)^*(g)$. We then composed the whole expression with $\Phi^t$ to have everything evaluated at the fixed point $x$. Now let $\{\p_i\}$ be an arbitrary (time independent!) basis  at $x$, then we have
\[
|\nabla u^t|^2 \circ \Phi^t (x) = |\nabla^{g^t} u|^2 \Big|_x = (g^t)^{ij} u_i u_j
\]
We can then compute using equation \eqref{MetricDerivatives}
\begin{align} \label{GLFirstVarComp}
\frac{d}{dt} E_{GL, \eps}(u^t) \Big|_{t = 0} &= \int_M \left[-  \frac{1}{2} \langle \dot{g}, S_{GL, \eps} \rangle + e_{GL, \eps}(u) \div(X)\right] dVol \\ \nonumber
&= \int_M \left[ - \langle S_{GL, \eps}, \nabla X \rangle + e_{GL, \eps}(u) \langle g, \nabla X \rangle \right] \\ \nonumber
&= \int_M \langle T_{GL, \eps}, \nabla X \rangle
\end{align}
having used symmetry of $S_{GL,\eps}$. This vanishes exactly when $u_{\eps}$ is a solution to \eqref{GLEquation}, reflecting the fact that the stress energy tensor is divergence free. \nl \nl
In the Abelian Yang--Mills--Higgs setting, define
\[
u^t = (\Phi^{-t})^*(u), \qquad \n^t = (\Phi^{-t})^*(\n)
\]
so that
\begin{align*}
E_{AYMH, \eps}(u^t, A^t) &= \int_M \left[|\nabla^t u^t|^2 + \eps^2 |F_{\n^t}|^2 + \frac{W(u^t)}{\eps^2} \right] d Vol \\
&= \int_M \left[ |\nabla^t u^t|^2 \circ \Phi^t + \eps^2 |F_{\n^t}|^2 \circ \Phi^t + \frac{W(u)}{\eps^2} \right] (\Phi^t)^*(d Vol)
\end{align*}
The analogous computation to equation \eqref{GLFirstVarComp} gives 
\begin{align*}
\frac{d}{dt} \left(|\nabla^t u^t|^2 \circ \Phi^t(x) \right) &= - \langle \dot{g}, \nabla u^* \nabla u \rangle \\
\frac{d}{dt} \left(|F_{\n^t}|^2 \circ \Phi^t(x) \right) &= - \langle \dot{g}, F_\n^* F_\n \rangle \\
\frac{d}{dt} E_{AYMH, \eps}(u^t, \n^t) \Big|_{t = 0}&= \int_M \eAYe(u,\n) \div(X) - \langle \dot{g}, S_{AYMH, \eps} \rangle \\ 
&= \int_M \langle T_{AYMH, \eps}, \nabla X \rangle
\end{align*}
\section{Computation of second inner variation}
In the Ginzburg--Landau case, we compute using \eqref{MetricDerivatives}:
\begin{align*}
\frac{d^2}{dt^2}\left( |\nabla u^t|^2 \circ \Phi^t \right)(x) \Big|_{t = 0}& = \ddot{g}^{ij} u_i u_j \Big|_x \\
&= \langle 2 (\dot{g} \circ \dot{g}) - \ddot{g}, \nabla u^* \nabla u \rangle 
\end{align*}
here,
\[
(\dot{g} \circ \dot{g})_{ij} := \dot{g}_{ik} g^{k \ell} \dot{g}_{\ell j}
\]
So that
\begin{align*}
\frac{d^2}{dt^2} \EGLe(u^t_{\eps}) \Big|_{t = 0} &= \int_M \frac{1}{2} \frac{d^2}{dt^2} \left( |\n^t u^t|^2 \circ \Phi^t \right) d Vol + \frac{d}{dt} \left( |\n^t u^t|^2 \circ \Phi^t \right) \Big|_{t = 0} \cdot \frac{d}{dt} (\Phi^t)^*(dVol_g) \Big|_{t = 0} + e_{\eps}(u) \frac{d^2}{dt^2} (\Phi^t)^*(dVol_g) \Big|_{t = 0} \\
&= \int_M \frac{1}{2}\langle 2 (\dot{g} \circ \dot{g}) - \ddot{g}, S_{GL, \eps} \rangle d Vol_M \\
& + \int_M - \langle \dot{g}, S_{GL, \eps} \rangle \div(X) d Vol \\
& + \int_M \eGLe(u_{\eps}) \left[ \div(\n_X X) - \Ric(X,X) - T_{TM}(X, X) + \div(X)^2\right] d Vol
\end{align*}
In the Abelian Yang--Mills--Higgs case, we proceed analogously:
\begin{align*}
\frac{d^2}{dt^2} \EAYe(u, A) \Big|_{t = 0} &= \int_M \frac{1}{2}\langle 2 (\dot{g} \circ \dot{g}) - \ddot{g}, S_{AYMH, \eps} \rangle d Vol_M \\
& + \int_M - \frac{1}{2} \langle \dot{g}, S_{AYMH, \eps} \rangle \div(X) d Vol \\
& + \int_M \eAYe(u_{\eps}) \left[ \div(\n_X X) - \Ric(X,X) - T_{TM}(X, X) + \div(X)^2 \right]
\end{align*}
We now reduce the second inner variation, using our lemmas about the stress-energy tensor. \nl \nl
\noindent \textbf{Proof of theorem }\ref{SecondInnerVarThm}: \nl
\noindent From lemma \ref{GLStressEnergyTensorConvergence} and equation \eqref{GLConvergenceAssumption}, we have that 
\begin{align*}
\lim_{\eps \to 0} \frac{d^2}{dt^2} \EGLe(u^t) \Big|_{t = 0} & = \pi \sum_{i = 1}^q m_i \int_{\Gamma_i} \Big( \frac{1}{2} \left[ 2\tr_{N \Gamma_i}(\dot{g} \circ \dot{g}) - \tr_{N \Gamma_i}(\ddot{g}) \right] - \tr_{N\Gamma_i}(\dot{g}) \div(X) \\
&  \qquad + [\div(\n_X X) - \Ric(X,X) - T_{TM}(X,X) + \div_{TM}(X)^2] \Big) d Vol
\end{align*} 
For $P, W  \subseteq TM$ subbundles and $X \in \Gamma(TM)$, let
\[
||\nabla_P^W X||^2 := \tr_P \langle \Pi^{W}\nabla_{\cdot} X, \Pi^{W} \n_{\cdot} X \rangle
\]
At a point $x \in M$, if $\{v_i\}_{i = 1}^k$ is an orthonormal basis for $P_x$ and $\{w_j\}_{j = 1}^{\ell}$ orthonormal for $W_x$, the above becomes
\begin{align*}
||\nabla_P^W X||^2\Big|_x &= \sum_{i = 1}^d \langle \Pi^W \n_{v_i} X, \Pi^W \n_{v_i} X \rangle \Big|_x \\
&= \sum_{i = 1}^d \sum_{j = 1}^{\ell} \langle \n_{v_i} X, w_j \rangle^2 \Big|_x 
\end{align*}
We compute in an orthonormal basis:
\begin{align*}
\tr_{N \Gamma_i}(\dot{g} \circ \dot{g}) &= \sum_{i = n-1}^n (\dot{g} \circ \dot{g})_{ii} \\
&= \sum_{i = n-1}^n \sum_{j = 1}^n \dot{g}_{ij}^2 \\
&= \sum_{i = n-1}^n \sum_{j = 1}^n [\langle \n_i X, \p_j\rangle^2 + \langle \n_j X, \p_i\rangle^2 + 2 \langle \n_i X, \p_j\rangle \langle \n_j X, \p_i\rangle ] \\
&= ||\nabla_{N \Gamma_i}^{TM} X||^2 + ||\nabla_{TM}^{N \Gamma_i} X||^2 +  2 \tr_{N \Gamma_i} \langle \n_{\n_{\cdot} X} X, \cdot \rangle\\
\tr_{N \Gamma_i}(\ddot{g}) &= \sum_{i = n-1}^n 2\langle \n_i \n_X X, \p_i \rangle + 2 \langle \n_i X, \n_i X \rangle - 2 R(X, \p_i, X, \p_i) \\
&= 2\div_{N \Gamma_i} (\n_X X) + 2 ||\nabla_{N \Gamma_i}^{TM} X ||^2 - 2 \tr_{N \Gamma_i}R(X, \cdot, X, \cdot) \\
\tr_{N \Gamma_i}(\dot{g}) &= 2\sum_{i = n-1}^n \langle \n_i X, \p_i \rangle \\
&= 2 \div_{N \Gamma_i} (X) 
\end{align*}
%
%
In sum
%
\begin{align*}
\lim_{\eps \to 0} \frac{d^2}{dt^2} \EGLe(u^t) \Big|_{t = 0} &= \pi \sum_{i = 1}^q m_i \int_{\Gamma_i} \Big(\left[ ||\nabla_{TM}^{N \Gamma_i} X||^2+ 2 \tr_{N \Gamma_i} \langle \n_{\n_{\cdot} X} X, \cdot \rangle - \div_{N \Gamma_i} (\n_X X) + \tr_{N \Gamma_i}R(X, \cdot, X, \cdot) \right] \\
& \qquad \quad - 2 \div_{N \Gamma_i}(X) \div(X) \\
& \qquad \quad + [\div(\n_X X) - \Ric(X,X) - T_{TM}(X,X) + \div_{TM}(X)^2] \Big) dVol_{\Gamma_i}
\end{align*}
With some cancellation, we group
\begin{align*}
\tr_{N \Gamma_i}R(X, \cdot, X, \cdot)  -\Ric(X, X) &= -\tr_{T \Gamma_i}R(X, \cdot, X, \cdot) \\
\div_{TM}(X)^2 - 2 \div_{N \Gamma_i}(X) \div(X)&= \div_{T \Gamma_i}(X)^2 - \div_{N \Gamma_i}(X)^2 \\
\div(\n_X X) - \div_{N \Gamma_i}(\n_X X) &= \div_{T \Gamma_i}(\n_X X)  \\
||\nabla_{TM}^{N \Gamma_i} X||^2 &= ||\nabla_{T \Gamma_i}^{N \Gamma_i} X||^2 + ||\nabla_{N \Gamma_i}^{N \Gamma_i} X||^2 \\
||\nabla_{N \Gamma_i}^{N \Gamma_i} X||^2 + 2 \tr_{N \Gamma_i}\langle \n_{\n_{\cdot} X} X, \cdot \rangle - T_{TM}(X,X) &= - T_{T \Gamma_i}(X,X) + \frac{1}{2} ||\dot{g}||_{N \Gamma_i}^2
\end{align*}
For $T$ as in equation \eqref{TTensorEqn}. This gives
\begin{align*}
\lim_{\eps \to 0} \frac{d^2}{dt^2} \EGLe(u^t) \Big|_{t = 0} &= \pi \sum_{i = 1}^q m_i \int_{\Gamma_i} \Big( - \tr_{T \Gamma_i}R(X, \cdot, X, \cdot) + \div_{T \Gamma_i}(X)^2 + \div_{T \Gamma_i}(\n_X X) + ||\nabla_{T \Gamma_i}^{N \Gamma_i} X||^2 - T_{T \Gamma_i}(X, X) \Big) \\
& \quad + \Big( - \div_{N \Gamma_i}(X)^2 + \frac{1}{2} ||\dot{g}||_{N \Gamma_i}^2 \Big) 	\\
&= \pi \sum_i m_i \left[ D^2A \Big|_{\Gamma_i}(X, X) + \int_{\Gamma_i} \Big( - \div_{N \Gamma_i}(X)^2 + \frac{1}{2} ||\dot{g}||_{N \Gamma_i}^2  \Big) d Vol  \right]
\end{align*}
having used equation \eqref{SecondVarArea}. Note that the error term is non-negative: let $\{\p_1,\p_2\}$ a basis for $N \Gamma_i$ which is orthonormal when restricted to $\Gamma_i$. We have:
\begin{align*}
- \div_{N \Gamma_i}(X)^2 + \frac{1}{2} ||\dot{g}||_{N \Gamma_i}^2 &= \sum_{i,j = n-1}^n - \langle \n_i X, \p_i \rangle \langle \n_j X, \p_j \rangle + \frac{1}{2}[\langle \n_i X, \p_j \rangle + \langle \n_j X, \p_i \rangle]^2 \\
&= [\langle \n_{n-1} X, \p_{n-1} \rangle - \langle \n_n X, \p_n \rangle]^2 + [\langle \n_{n-1} X, \p_n \rangle + \langle \n_n X, \p_{n-1} \rangle]^2 \\
&= \frac{1}{4} [\dot{g}_{n-1, n-1} - \dot{g}_{n,n}]^2 + \dot{g}_{n-1,n}^2 \\
&= [\langle \n_{n-1} X, \p_{n-1} \rangle - \langle \n_n X, \p_n \rangle]^2 + [\langle \n_{n-1} X, \p_n \rangle + \langle \n_n X, \p_{n-1} \rangle]^2
& \geq 0
\end{align*}
\noindent We note that the error term vanishes when $\nabla_i X \Big|_{\gamma} = 0$, e.g. if $X \Big|_{\Gamma_i} = f^1(s) \p_1 + f^2(s) \p_2$, where $s$ is a coordinate on $\Gamma_i$, and $\p_1, \p_2$ are normal vector fields corresponding to a parallel frame on $\Gamma_i$. Such a vector field can be extended to a vector field on $M$ via a bump function in a tubular neighborhood of $\Gamma_i$ (see \cite{gaspar2020second}, Appendix). We also note that this is the same error term in Cheng \cite{ChengInstability}, Proposition 2.6, as well as Le \cite{le2015second}, Theorem 1.5. \nl \nl
The reduction of $\lim_{\eps \to 0} \frac{d^2}{dt^2} E_{AYMH, \eps}(u^t, \n^t)$ is identical with the same result (but a normalizing factor of $2\pi$) since the reduction only depends on the limit of the energy measure, the stress energy tensor, and derivatives of $g^t$.
\section{Proof of Theorem \ref{MorseIndexTheorem}}
This section takes strong inspiration from theorem $A$ in \cite{gaspar2020second}, \S 4, and Lemma 3.12 \cite{hiesmayr2018spectrum}. \nl \nl
\noindent Recall assumption \ref{ConvergenceAssumption} so that 
\[
\mu_V = \sum_i m_i d \mathcal{H}^{n-2}_{\Gamma_i}
\]
is the measure associated to our limiting varifold, $V$. As in \cite{hiesmayr2018spectrum} (main theorem) and \cite{gaspar2020second} (theorem A), fix a set $U \subset \subset M \backslash \text{sing}(V)$. Define 
\[
Q_V(X) = D^2 V(X) = \sum_{i = 1}^q m_i \int_{\Gamma_i} |\nabla^{\perp} X|^2 - \Ric(X, X) - |\langle A (\cdot, \cdot), X \rangle|^2
\]
for $X = \{X_i\}$ a collection of normal vector fields to each $\Gamma_i$.
Given a sequence of $\{u_{\eps_k}\}$ solutions to equation \eqref{GLEquation} and $\{(u_{\eps_k}, \n_{\eps_k})\}$ solutions to equations \eqref{AYMHEquation} satisfying our assumptions \ref{ConvergenceAssumption}, we have 
\begin{align*}
\lim_{k \to \infty} \mu_{GL, \eps_k} &= \mu_V \\
\lim_{k \to \infty} \mu_{AYMH, \eps_k} &= \mu_V
\end{align*}
in a weak sense. Let $\tilde{X}$ be an extension of $X$ such that $\nabla_{N \Gamma_i} \tilde{X} \Big|_{\Gamma_i} = 0$ (see, e.g. \cite{gaspar2020second}, Appendix for details). Let $\{u^t_{\eps_k}\}$ (resp. $\{(u^t_{\eps_k}, \n^t_{\eps_k})\}$) be the corresponding family of functions (resp. functions and connections) after pulling back our sequence of Ginzburg--Landau (resp. Abelian Yang--Mills--Higgs) solutions by $\Phi^t$, the flow associated to $\tilde{X}$ of on $M$. Then we have
\[
\frac{1}{2\pi} \lim_{k \to \infty} \frac{d^2}{dt^2} \EAYe(u^t_{\eps_k}, \n^t_{\eps_k}) \Big|_{t = 0} = \frac{1}{\pi}\lim_{k \to \infty} \frac{d^2}{dt^2} \EGLe(u^t_{\eps_k}) \Big|_{t = 0} = Q_V(X)
\]
Denote
\begin{align*}
Q_{\eps, GL}(\tilde{X}) &:= \frac{d^2}{dt^2} \EGLe(u^t_{\eps_k}) \Big|_{t = 0} = \EGLe''(u_{\eps_k})(\langle \nabla u_{\eps_k}, \tilde{X} \rangle, \langle \nabla u_{\eps_k}, \tilde{X} \rangle) \\
Q_{\eps, AYMH}(\tilde{X}) &:= \frac{d^2}{dt^2} \EAYe(u^t_{\eps_k}, \n^t_{\eps_k}) \Big|_{t = 0} = \EAYe''(u_{\eps_k}, \n_{\eps_k})(\nabla_{\eps_k,\tX} u_{\eps_k}, F_{\n_{\eps_k}}(\tX, \cdot)) 
\end{align*}
%
%
We now recall the variational definition of $\lambda_p^{\eps}(U)$, of the second variation operators for Ginzburg--Landau and Abelian Yang--Mill--Higgs.
\begin{align} \nonumber
\lambda_{GL, p}^{\eps}(U) &= \inf_{\dim E = p} \sup_{\phi \in E \backslash 0} \frac{\EGLe(u_{\eps})''(\phi, \phi)}{||\phi||_{L^2(M)}^2} \\ \nonumber
\lambda_{AYMH, p}^{\eps}(U) &= \inf_{\dim E = p} \sup_{(\phi,a) \in E \backslash 0} \frac{\EAYe(u_{\eps}, \n_{\eps})''((\phi,a), (\phi,a))}{||(\phi,a)||_{L^2(M)}^2} 
\end{align}
In $\lambda_{GL, p}$, $E$ consists of functions $\phi: M \to \C$, whereas in $\lambda_{AYMH, p}$, $E$ consists of sections, $\phi \in \Gamma(L)$, and $i\R$ valued one forms, $a \in \Omega^1(TM, i \R)$. We similarly define $\lambda_p(U)$ for the second variation of area
\begin{align} \nonumber
\lambda_p(U) &= \inf_{\dim E = p} \sup_{X \in E \backslash 0} \frac{\sum_{i  =1}^q \int_{\Gamma_j} ||\nabla^{\perp}X||^2 - \tr_{T \Gamma_i}(R(X, \cdot, X, \cdot)) - |\langle A(\cdot, \cdot), X \rangle|^2}{||X||_{L^2(U \cap \text{reg}(V))}^2} \\ \label{WeightedCharacterization}
&= \inf_{\dim E = p} \sup_{X \in E \backslash 0} \frac{Q_V(X)}{||X||_{L^2(V)}^2}
\end{align}
for 
\[
||X||_{L^2(V)}^2 = \sum_{i  =1}^q m_i \int_{\Gamma_i} |X|^2 d \mathcal{H}^{n-2}_{\Gamma_i}
\]
Note that equation \eqref{WeightedCharacterization} holds by normalizing each of $X\Big|_{\Gamma_i}$ by $\sqrt{m_i}$, as noted by Hiesmayr \cite{hiesmayr2018spectrum}, \S 3.2. \nl \nl
\noindent \textbf{Proof of theorem } \ref{MorseIndexTheorem}: \nl
Given $\delta > 0$, there is $E$, a $p$-dimensional space spanned by $X_1, \dots, X_p$, and $\{c_i\}_{i = 1}^p$ so that 
\begin{equation} \label{LambdaJacobiLowerBound}
\max_{\{c_i\} \in S^{p - 1}}\frac{Q_V(\sum_{i} c_i X_i)}{||\sum_i c_i X_i||^2_{L^2(V)}} \leq \lambda_{p}(W) + \delta
\end{equation}
Let $\{\tX_i\}$ be normal, parallel extensions of each $X_i$ to vector fields on all of $TU$ with compact support and $\nabla^{N \Gamma_i} \tX_j \Big|_{\Gamma_i} = 0$ for all $i,j$. For any $c = \{c_i\} \in S^{p-1}$ we will denote 
\[
c \cdot \tilde{X} := \sum_{i = 1}^p c_i \tX_i
\]
\noindent For arbitrary $b = \{b_i\}$, consider the maps
\begin{align*}
F_{GL,\eps_k}&: S^{p-1} \to H^1_0(U) \\
b \cdot \tilde{X} &\mapsto \langle b \cdot \tilde{X}, \nabla u_{\eps_k} \rangle \\
F_{AYMH,\eps_k}&: S^{p-1} \to H^1_0(U, L) \times H^1_0(\Omega^1(TU, L)) \\
b \cdot \tilde{X} &\mapsto (\nabla_{\eps_k, b \cdot \tilde{X}} u_{\eps_k}, F_{\n_{\eps_k}}(b \cdot \tilde{X}, \cdot))
\end{align*}
%
Note that for $\eps$ sufficiently small (or $k$ sufficiently large), the maps $F_{GL, \eps_k}$ and $F_{AYMH, \eps_k}$ are injective. To see this, suppose not, then there exists a sequence of $\{c^k = (c_1^k, \dots, c_p^k)\}$ such that 
\[
F_{GL, \eps_k}(c^k \cdot \tilde{X}) = 0
\]
Since $\{c^k \} \subseteq S^{p-1}$, then up to relabelling it with a subsequence, we have $c^{k} \to c$. However, the weak convergence of the stress energy tensor in equation \eqref{AYMHStressTensorEq} gives
\begin{align*}
0 &= \lim_{k \to \infty} \langle F_{GL,\eps_{k}}(c^{k} \cdot \tilde{X}) , F_{GL, \eps_{k}}(c^k \cdot \tilde{X}) \rangle \\
&= \lim_{k \to \infty} \langle (c^k \cdot \tilde{X}) \otimes (c^k \cdot \tilde{X}), S_{GL, \eps_k} \rangle \\
&= \sum_{j =1}^q m_j \int_{\Gamma_j} \tr_{N \Gamma_i}((c \cdot \tilde{X}) \otimes (c \cdot \tilde{X})) \\
&= \sum_{j =1}^q m_j \int_{\Gamma_j} ||\Pi^{N \Gamma_i} (c \cdot \tilde{X})||^2
\end{align*}
But this immediately implies $c \cdot \tilde{X} = 0$, a contradiction because $c \in S^{p-1}$. The same computation works for Abelian Yang--Mills--Higgs: suppose we have a sequence $\{c^k\}$ such that 
\[
F_{AYMH, \eps_k}(c^k \cdot \tilde{X}) = 0 
\] 
Using the inner product on pairs of sections and one forms \eqref{EpsInnerProduct} and the convergence in equation \eqref{GLStressTensorEq}, we have
\begin{align*}
0 &= \lim_{k \to \infty} \langle F_{AYMH,\eps_{k}}(c^{k} \cdot \tilde{X}) , F_{AYMH, \eps_{k}}(c^k \cdot \tilde{X}) \rangle \\
&= \frac{1}{2}  \lim_{k \to \infty} \langle (c \cdot \tilde{X}) \otimes (c \cdot \tilde{X}), S_{AYMH, \eps_k} \rangle \\
&= \frac{1}{2} \sum_{j =1}^q m_j \int_{\Gamma_j} \tr_{N \Gamma_i}((c \cdot \tilde{X}) \otimes (c \cdot \tilde{X})) \\
&= \frac{1}{2} \sum_{j =1}^q m_j \int_{\Gamma_j} ||\Pi^{N \Gamma_i} (c \cdot \tilde{X})||^2
\end{align*}
Thus, the spaces of 
\begin{align*}
W_{GL, \eps_k, p} &= \{\langle b \cdot \tilde{X}, \nabla u_{\eps_k} \rangle \; | \; b \in S^{p-1} \} \\ 
W_{AYMH, \eps_k, p} &= \{ (\nabla_{\eps_k, b \cdot \tilde{X}} u_{\eps_k}, F_{\eps_k}(b \cdot \tilde{X}, \cdot)) \; | \; b \in S^{p-1}\}
\end{align*}
are both $p$-dimensional linear subspaces of $H_0^1(U)$. By definition, we then have
\begin{align} \label{GLLambdaBound}
\lambda_{GL, p}^{\eps_k}(U) & \leq \sup_{c \in S^{p - 1}} \frac{Q_{GL, \eps_k}(c \cdot \tilde{X})}{||\langle c \cdot \tilde{X}, \nabla u_{\eps_k} \rangle||_{L^2(M)}^2} \\ \label{AYMHLambdaBound}
\lambda_{AYMH, p}^{\eps_k}(U) & \leq \sup_{b \in S^{p - 1}} \frac{Q_{AYMH, \eps_k}(b \cdot \tilde{X})}{||(\nabla_{\eps_k,c \cdot \tilde{X}} u_{\eps_k}, F_{\n_{\eps_k}}(b \cdot \tilde{X}, \cdot)) ||_{L^2(M)}^2} 
\end{align}
We define
\begin{align*}
\mu_{GL, p}(U) &:= \limsup_{k \to \infty} \lambda_{GL, p}^{\eps_k}(U) \\
\mu_{AYMH, p}(U) &:= \limsup_{k \to \infty} \lambda_{AYMH, p}^{\eps_k}(U)
\end{align*}
Up to taking a subsequence, we can replace $\limsup$ with $\lim$ in the above. Now for each $k$, choose $c^k$ and $b^k$ which achieves the $\sup$ in equations \eqref{GLLambdaBound} and \eqref{AYMHLambdaBound}. This means for our choice of $\delta > 0$, there exists $K$ such that for $k > K$ there exists $\{c^k\}$ and $\{b^k\}$ so that
\begin{align*}
\mu_{GL,p}(U)  - \delta &\leq \frac{Q_{GL, \eps_k}(c^k \cdot \tX)}{||\langle c^k \cdot \tX, \nabla u_{\eps_k} \rangle||_{L^2(M)}^2} \\
\mu_{AYMH,p}(U)  - \delta &\leq \frac{Q_{AYMH, \eps_k}(b^k \cdot \tX)}{||(\nabla_{\eps_k,b^k \cdot \tX} u_{\eps_k}, F_{\n_{\eps_k}}(b^k \cdot \tX, \cdot)) ||_{L^2(M)}^2}
\end{align*}
\noindent Now taking another subsequence $c^k \to c$, $b^k \to b$ and using the convergence of the second inner variation in equation \ref{SecondInnerVarEqn}, we have 
\begin{align*}
\lim_{k \to \infty} Q_{GL, \eps_k}(c^k \cdot \tX) &= \pi Q_V(c \cdot \tilde{X}) \\
\lim_{k \to \infty} Q_{AYMH, \eps_k}(b^k \cdot \tX) &= 2\pi Q_V(b \cdot \tilde{X}) 	
\end{align*}
But by convergence of stress tensors (\ref{AYMHStressEnergyTensorConvergence}, \ref{GLStressEnergyTensorConvergence}), we also have
\begin{align*}
\lim_{k \to \infty} ||\langle c^k \cdot \tX, \nabla u_{\eps_k} \rangle||_{L^2(M)}^2 & = \lim_{k \to \infty} \langle (c^k \cdot \tX) \otimes (c^k \cdot \tX), S_{GL, \eps_k} \rangle \\ 
&= \pi \sum_{j = 1}^{q} m_j \int_{\Gamma_j}  \tr_{N \Gamma_j}\left(c \cdot \tX \right)^2 \\
&= \pi \sum_{j = 1}^{q} m_j \int_{\Gamma_j}  ||c \cdot \tX||^2 \\
&= \pi ||c \cdot \tX||_{L^2(V)}^2 \\
\lim_{k \to \infty} ||(\nabla_{\eps_k,b^k \cdot \tX} u_{\eps_k}, F_{\n_{\eps_k}}(b^k \cdot \tX, \cdot)) ||_{L^2(M)}^2 & = \frac{1}{2} \lim_{k \to \infty} \langle (c^k \cdot \tX) \otimes (c^k \cdot \tX), S_{AYMH, \eps_k} \rangle \\  
& = \frac{1}{2} 2 \pi ||b \cdot \tX||_{L^2(V)}^2 \\
& = \pi ||b \cdot \tX||_{L^2(V)}^2 
\end{align*}
%
This tells us that for $k$ sufficiently large and using equation \eqref{LambdaJacobiLowerBound}:
\begin{align*}
\mu_{GL, p}(U) - \delta & \leq \frac{Q_V(c \cdot \tX)}{||c \cdot \tX||_{L^2(V)}^2} \leq \lambda_{p}(U) + \delta \\
\mu_{AYMH, p}(U) - \delta  & \leq \frac{2 Q_V(b \cdot \tX)}{||b \cdot \tX||_{L^2(V)}^2} \leq 2 \lambda_{p}(U) + \delta
\end{align*}
We note that the factor of $2$ is an artifact of the convention of equation \eqref{AYMHEnergy}, i.e. it would not appear if we considered half of the Abelian Yang--Mills--Higgs energy instead. Since $\delta$ is arbitrary, we conclude 
\begin{align*}
\mu_{GL, p}(U) &\leq \lambda_p(U) \\
\mu_{AYMH, p}(U) &\leq 2\lambda_p(U)
\end{align*}
Recalling that 
\begin{align*}
\text{Ind}(\text{reg}(V))&= \sup_U \# \{p \in \mathcal{N} \; | \; \lambda_p(U) < 0\} \\
\implies \text{Ind}(\text{reg}(V)) & \leq \limsup_{k \to \infty} \text{Ind}_{GL}(u_{\eps_k}) \leq M \\
\implies \text{Ind}(\text{reg}(V)) & \leq \limsup_{k \to \infty} \text{Ind}_{AYMH}(u_{\eps_k}, \n_{\eps_k}) \leq M
\end{align*}
We see that the morse index of the varifold is bounded above by the upper bound on the morse index of our sequence of solutions as $\eps_k \to 0$. 
%
%
\bibliography{Second_Inner_Variation}{}
\bibliographystyle{plain}
\end{document}